\documentstyle[12pt,amssymbols]{article}  
\newcommand{\bq}{\begin{equation}}  
\newcommand{\eq}{\end{equation}}  
\newcommand{\bqa}{\begin{eqnarray}}  
\newcommand{\eqa}{\end{eqnarray}}  
\newcommand{\ra}{\rightarrow}

\def\half{{1 \over 2}}  
\def\s{\sigma}

\def\ep{\epsilon}

\def\ov{\over}

\def\ed{\end{document}}  
  
\def\ra{\rightarrow}  
  
\def\2pi{1\over 2\pi i}  

\def\~{\tilde}
\def\newline{\hfil\break}

\def\ra{\rightarrow}

\def\sq2{\sqrt{2}}  
\def\sqk2{\sqrt{2(k+2}}  
\def\sqk{\sqrt{k}}

\def\be{\begin{equation}}  
\def\ee{\end{equation}}  
\def\br{\begin{array}}  
\def\er{\end{array}}  
\def\bea{\begin{eqnarray}}  
\def\eea{\end{eqnarray}}  
\def\ba{\begin{equation}\begin{array}}  
\def\ea{\end{array}\end{equation}}  
\def\bac{\begin{equation}\begin{array}{rll}}

\newcommand{\uq}{U_q (\widehat{sl(2)})}

\def\Z{{\Bbb Z}}

\def\pl{\prod\limits}

\def\ep{\epsilon}

\begin{document}  
\rightline{ITP-SB-96-21}  
\rightline{ April, 1996}  
\vbox{\vspace{-10mm}}  
\vspace{1.0truecm}  
\begin{center}  
{\LARGE \bf   
Exact two-spinon dynamical 
correlation function
of the Heisenberg model
 }\\[8mm]  
{\large A.H. Bougourzi$^{1}$, M. Couture$^{2}$ and M.
Kacir$^{3}$  }\\  
[6mm]{\it 
$^{1}$Institute of Theoretical Physics\\  
SUNY at Stony Brook\\  
Stony Brook, NY 11794
\\[3mm]

$^{2}$Neutron \& Condensed Matter Sciences\\
AECL Research, Chalk River Laboratories\\
Chalk River, Ontario, Canada KOJ 1J0\\[3mm]
$^{3}$Service de Physique Theorique\\
 CE-Saclay\\
F-91191 Gif-sur-Yvette Cedex, France
}\\[20mm]  
  
\end{center}  
\vspace{1.0truecm}  
\begin{abstract}  
We derive the exact contribution of   two spinons to the 
dynamical correlation
function of the spin-1/2  Heisenberg model.
For this, we use the isotropic limits of the exact form
factors that have been recently computed  through the
quantum affine symmetry of the anisotropic Heisenberg 
model $XXZ$.   
\end{abstract}  
\newpage  
\section{Introduction}

Ever since Niemeijer derived an exact expression for the
dynamical correlation function (DCF) at any temperature 
 of the spin-1/2 $XY$
model \cite{Nie67}, there has been a considerable amount of 
work in trying
to extend his result to the more physically interesting case
of the  isotropic spin-1/2 Heisenberg model (i.e., $XXX$ model).
For more details on the definitions of the spin chain models see
Refs. \cite{Bax82,Kor93}. 
However, so far only approximate but very accurate results have
been computed. For a
comprehensive 
  historic review, the importance of the DCF, 
an account of the existing results, 
 and a list of references 
on this subject we recommend Refs. 
\cite{Rolal86,Mulal81,MuSh84,Aff89b,Ten95,Tenal95}. 
In particular, 
in Ref. \cite{Mulal81} an ansatz for
the DCF at zero temperature was proposed based on 
Niemeijer's result and other 
approximate numerical 
and analytical results. But it has never been established
whether this ansatz includes contributions  from just 
two spinons
or more. The spinon picture has been rigorously 
studied in the case of the Heisenberg model 
in \cite{FaTa79}, and is of great theoretical and 
experimental interest. 
Let us just mention that the main stumbling block in trying to
compute the eaxct DCF  for the Heisenberg model is due to the
absence of exact results for the form factors. Unfortunately,
so far the powerful method of the Bethe ansatz has not been able
to provide them. The form factors are well understood
now just in two-dimensional quantum field theories with familiar
relativistic dispersion relations \cite{Smi92} 
but not yet in lattice 
models.  However,  an approach based on the
concept of exact resolution of
dynamics through just
infinite-dimensional symmetries, and 
which is widely used
in the context of string theory and conformal field theory, has
been recently upgraded to  
 the ``massive" $XXZ$ model
in  \cite{Daval93,JiMi94}. 
It provides almost all exact physical quantities (static
correlation functions and form factors). 
Its only shortcomings are that it is based on a
relatively complicated
symmetry, which is the quantum affine algebra $\uq$, 
and that the
 physical quantities  it leads to
are somewhat complicated to deal with. To be more precise,
they typically
have a multi-integral form. Our main point 
in this paper is that all these
latter problems disappear in the particular case of the 
two-spinon
form factors of the Heisenberg model because they 
have a simpler form.
Therefore we use them to compute the more interesting quantity 
of the exact two-spinon DCF at zero temperature.

Our letter is organized as follows:
first we briefly review the results related to the
diagonalization of the anisotropic $XXZ$ Heisenberg model
 following Ref. \cite{JiMi94}. Then we define the two-spinon
DCF in the case of $XXZ$ in terms of the form factors of this
model. Finally, we take the isotropic limit. 
As mentioned earlier,
the crucial 
point is that considerable simplifications take place 
in this case
due to the isotropy of the Heisenberg model, and thus 
 allow us
to derive a simple  exact  formula for DCF of
two-spinon. We hope that the exact results contained in this
paper will shed some new light on the spinon picture from both
theoretical and experimental perspectives where already a
considerable amount of work  has been achieved over the past 
few decades.

\section {Diagonalization of 
the anisotropic Heisenberg model}

The Hamiltonian of the anisotropic (XXZ) Heisenberg model 
is defined by
\be { H_{XXZ}}=
-\half \sum_{n=-\infty}^{\infty} (\s_n^x \s_{n+1}^x 
+\s_n^y \s_{n+1}^y + \Delta\s_n^z \s_{n+1}^z) 
\label{hamiltonian}, 
\ee 
where $\Delta=(q+q^{-1})/2$ is the anisotropy parameter. 
Here $\sigma_n^{x,y,z}$ 
are the usual
Pauli matrices acting at the $n^{\rm th}$ position of the 
formal 
infinite tensor product 
\be W= \cdots V \otimes V \otimes V  \cdots 
\label{infprod},\eq 
where $V$ is the two-dimensional representation of $U_q(sl(2))$ 
quantum group. 
We consider the model in 
the anti-ferromagnetic regime $\Delta <-1$, i.e., 
$-1<q<0$.
The  main point of Refs. \cite{Daval93,JiMi94}  
is that the action of $H_{XXZ}$ 
on $W$ is
not well defined due to the appearance of 
divergences. However, this model is symmetric 
under the quantum group $\uq$, and therefore the eigenspace
is identified with the following level 0 $\uq$ 
module:
\be
{\cal F}=\sum_{i,j}V(\Lambda_i)\otimes V(\Lambda_j)^*, 
\ee
where $\Lambda_i$ and $V(\Lambda_i); i=0,1$ are  level 
1 $\uq$-highest
weights and $\uq$-highest weight modules, respectively.
Roughly speaking, $V(\Lambda_i)$ is identified 
with the
subspace of the formal semi-infinite space
\be X= \cdots V \otimes V \otimes V,\eq 
consisting  of all linear combinations of  
spin configurations with  fixed boundary conditions such that 
the eigenvalues of 
$\sigma^z_n$ are $(-1)^{i+n}$ in the limit $
 n\ra -\infty$. The 
particle picture 
of this Hamiltonian is given in terms of vertex operators which
act as intertwiners of $\uq$ modules, and which create 
the set of eigenstates (spinons) $\{|\xi_1,\cdots 
\xi_n>_{\ep_1,\cdots \ep_n;i}, n\geq 0\}$. Here $i$ fixes
the boundary conditions, $\xi_j$ are the spectral parameters
living on the unit circle,
and $\ep_j=\pm 1$ are 
the spins of the
spinons. The actions of $H_{XXZ}$ and the translation operator 
$T$, which shifts the spin  chain by one site, on ${\cal F}$ 
are given by 
\bac T|\xi_1,\cdots,\xi_n>_i &=&\pl_{i=1}^n\tau(\xi_i)^{-1}
|\xi_1,\cdots,\xi_n>_{1-i},\quad 
T|0>_i =|0>_{1-i}, \\
H_{XXZ}|\xi_1,\cdots,\xi_n>_i 
&=&\sum_{i=1}^n e(\xi_i)|\xi_1,\cdots,\xi_n>_i, 
\label{states}\ea  
where
\bac \tau(\xi)&=& \xi^{-1} 
{\theta_{q^4} (q \xi^2) \ov \theta_{q^4} (q 
\xi^{-2})}=e^{-i p(\alpha)},\quad p(\alpha)=am({2K\over \pi}
\alpha)-{\pi/2}, \\
e(\xi) &=&{ 1-q^2 \ov 2 q} \xi {d \ov d \xi} 
\log \tau(\xi)= {2K\over \pi} \sinh({\pi K^\prime\over K})
dn({2K\over \pi}\alpha). \label{enmom}
\ea
Here, $e(\xi)$ and  $p(\alpha)$ are the energy and the 
momentum of
the spinon respectively,  $am(x)$ and $dn(x)$ are the usual 
elliptic amplitude and
delta functions, with the complete elliptic integrals
 $K$ and $K^\prime$, and 
\bac
q&=&-\exp(-\pi K^\prime/K),\\
\xi&=&ie^{i\alpha},\\ 
\theta_x(y)&=&(x;x)_{\infty} (y;x)_{\infty} 
(x y^{-1};x)_{\infty},\\
(y;x)_{\infty}&=&\prod_{n=0}^{\infty} (1-y x^n).
\ea
This means, $\sigma^{x,y,z}(t, n)$ at time $t$ and position
$n$ are related to $\sigma^{x,y,z}(0,0)$ at time 0 and position
0 through:
\be \sigma^{x,y,z} (t,n)=\exp(i tH_{XXZ} ) T^{-n} 
\sigma^{x,y,z} 
(0,0) T^{n} 
\exp(-i tH_{XXZ} ). \ee

The completeness relation reads \cite{JiMi94}: 
\be {\Bbb I}=\sum_{i=0,1}\sum_{n \geq 0} \sum_{\ep_1,\cdots,
\ep_n=\pm 1}
{1 \ov {n !}} \oint  {d\xi_1\over 2\pi i \xi_1} \cdots 
 {d\xi_n\over 2\pi i \xi_n}
 |\xi_n,\cdots,\xi_1>_{{\ep_n,\cdots,\ep_1};i}\;
{_{i;{\ep_1,\cdots,\ep_n}}{<\xi_1,\cdots,\xi_n|}}. 
\ee

\section{Two-spinon dynamical correlation function of the 
Heisenberg model}
First, we will define the dynamical correlation function we are
considering in the case of the anisotropic Heisenberg model,
where the particle picture is well understood and the 
form factors are known exactly  \cite{JiMi94}. However, let us
note that the expressions of these form factors are very
complicated to lead to a closed formula for the
DCF in the anisotropic case. 
But, in the isotropic limit,
i.e., $q\ra -1$
one of  them simplifies substantially, and using the fact 
that all the nonvanishing components of the 
DCF
are equal, we find the same closed formula for all of them in
this limit. 
 
Let us recall the definition of one of the components of
 the DCF in the case of the
XXZ
model. Up to an overall conventional 
normalization factor it is given by
\be
S^{i,+-}(w, k)=
\int_{-\infty}^{\infty} dt \sum_{n\in\Z}
e^{i(wt+kn)} {_i}< 0|\sigma^+(t, n)\sigma^-(0,0)|0>_i,
\ee
here $w$ and $k$ are the energy and momentum transfer
respectively, and $i$ corresponds to the boundary condition. 
Later we will
find that the DCF is in fact independent of $i$.
Using the completeness relation, 
the two-spinon contribution is given by
\bac
S_2^{i,+-}(w,k) &=& {\pi} 
\sum_{n\in \Z} \sum_{\ep_1,\ep_2} 
\oint {d\xi_1\over 2\pi i \xi_1}
 {d\xi_2\over 2\pi i \xi_2}
\exp\left(in(k+p(\xi_1)+p(\xi_2))\right)
 \delta (w-e(\xi_1)
-e(\xi_2))\\
&&\times {_{i+n}<0|}\sigma^+ 
(0,0)|\xi_2,\xi_1>_{\ep_2,\ep_1;i+n}\:
{_{i;\ep_1,\ep_2}<\xi_1,\xi_2|}
\sigma^-(0,0)|0>_i.
\ea
This can be re-written as
\bac
&& S_2^{i,+-}(w,k)=\pi
\sum_{\ep_1,\ep_2} 
 \oint {d \xi_1\over 2\pi i \xi_1}
{ d \xi_2\over 2 \pi i \xi_2} 
  \sum_{n\in \Z} \exp\left( 2 i n(k+p(\xi_1)+p(\xi_2))\right)
\delta (w-e(\xi_1)-e(\xi_2))\\
&&\times \left(
{_{i}<0|}\sigma^+(0,0)|\xi_2,\xi_1>_{\ep_2,\ep_1;i}
\:
{_{i;\ep_1,\ep_2}<\xi_1,\xi_2|}
\sigma^-{(0,0)}|0>_i \right.\\
&& +
\left. \exp\left( i(k+p(\xi_1)+p(\xi_2)\right)
{_{1-i}<0|}\sigma^+(0,0)|\xi_2,\xi_1>_{\ep_2,
\ep_1;1-i}\:
{_{i;\ep_1,\ep_2}<\xi_1,\xi_2|}
\sigma^-{(0,0)}|0>_i \right).
\label{corr2}\ea
The non-vanishing form factors have
been computed in \cite{JiMi94}, and satisfy the
following relations:
\bac {_i{<}}0|\sigma^-(0,0)|\xi_2,\xi_1>_{++;i} &=& 
{_{1-i}{<0|}}\sigma^+(0,0)|\xi_2,\xi_1>_{--;1-i}, \\
{_i{<}}0|\sigma^-(0,0)| \xi_2^{*},\xi_1^{*}>_{++;i}&=&
{_i{<}}0|\sigma^+(0,0)| -q \xi_1,-q \xi_2>_{--;i},\\
{_{i;--}{<}} \xi_1,\xi_2|\sigma^-(0,0)|0>_i&=&
{_i{<0|}}\sigma^-(0,0)|-q \xi_1,-q \xi_2>_{++;i}.
\label{FF}\ea

Now the isotropic limit $q\ra -1$ is performed
by first  making the following redefinitions:
\bac
\xi&=&ie^{{\epsilon\beta\over i\pi}},\\
q&=&-e^{-\epsilon},\quad\epsilon\ra 0^+.
\ea
Here $\beta$ is the appropriate 
spectral parameter  for the
Heisenberg model.

Then, one finds the following exact isotropic 
limits \cite{JiMi94}:
\bac
&&|{_{i}{<0|}}\sigma^+(0,0)|\xi_2,\xi_1>_{--;i}|^2
{d\xi_{1}\over 2\pi i  \xi_{1}} 
{d\xi_{2}\over  2\pi i \xi_{2}}\ra
{ \Gamma(3/4)^2|A_{-}(\beta_1-\beta_2)|^2\over
16\Gamma(1/4)^2 |A_+(i\pi/2)|^2 |A_{-}(i\pi/2)|^2
\cosh(\beta_1)\cosh(\beta_2)} d\beta_1 d\beta_2,\\
&&p(\xi)\ra p(\beta),\quad {\rm
s.t.}\quad \cot(p(\beta))=\sinh(\beta),\quad -\pi\leq p(\beta)
\leq 0,\\
&&e(\xi)\ra e(\beta)={\pi\over \cosh(\beta)}=-\pi \sin(p(\beta)),
\ea
where
\be
|A_{\pm}(\alpha)|^2=
\exp\left( -\int_{0}^{\infty}dx {(\cosh(2 x(1-{\delta\over\pi}))
\cos({2 x \gamma
\over \pi})-1)
\exp(\mp x)\over x \sinh(2 x)\cosh(x)}\right).
\ee
Here $\Gamma(x)$ is the usual gamma function and 
$\alpha=\gamma+i\delta$, with $\gamma$ and $\delta$ being
real.
Restricting
to the first Brillouin zone, integrating  the continuous and
discrete delta functions,
keeping track of the Jacobian factors, the
energy-momentum  conservation
relations, and the isotropic limits of (\ref{FF}), we find that 
$S^{i,+-}_2(w, k)$ is independent of $i$ (which is henceforth 
omitted ) and 
simplifies substantially to:
\be
S^{+-}_2(w, k)={\pi^2 {\Gamma(3/4)^2}
\Theta(2\pi \sin(k/2)-w)\Theta(
w-\pi|\sin(k)|)\over 4 \Gamma(1/4)^2
|A_{-}(i\pi/2)|^2|A_{+}(i\pi/2)|^2}
{|A_{-}(\bar\beta_1-\bar\beta_2)|^2\over
\sqrt{(2\pi \sin(k/2))^2-w^2}},
\label{exact}
\ee
where $\Theta$ is the Heaviside step function, and 
for fixed $w$ and $k$,  
$\bar\beta_1$ and $\bar\beta_2$ are the solutions to:
\bac
w&=&e(\bar\beta_1)+e(\bar\beta_2),\\
k&=&-p(\bar\beta_1)-p(\bar\beta_2).
\ea
Note that the pair $(\bar\beta_1, \bar\beta_2)$ 
is identified with the pair
$(\bar\beta_2, \bar\beta_1)$.

Let us now make some comments about $S^{+-}_2(w,k)$ as given 
by (\ref{exact}). From the isotropy of the Heisenberg model 
and the inclusion 
of both sectors $i=0$ and $i=1$,
we obtain all the non-vanishing components of the DCF from
$S^{+-}_2(w,k)$ as:
\be
S^{xx}_2(w,k)=S^{yy}_2(w,k)=S^{zz}_2(w,k)=4S^{+-}_2(w,k),
\ee 
with
\be
\sigma^{\pm}={\sigma^x\pm i\sigma^y\over 2}.
\ee 
Furthermore, from the dispersion relations of two spinons,  
$w$, as a function of $k$, lies between two boundaries: the lower
one is given by the famous des Cloizeaux-Pearson dispersion 
relation, i.e.,
\be
w_l=\pi|\sin(k)|,
\ee
whereas the upper one is given by the dispersion relation
\be
w_u=2\pi\sin(k/2).
\ee

Note that despite its square root singularity, $S_2^{+-}(w,k)$
actually vanishes in the vicinity of the upper boundary.
Moreover, it diverges in the vicinity of the lower boundary.
It would be very interesting to compare our results
with presently existing approximate results, and especially
the  ansatz made for the two-spinon DCF 
in Ref. \cite{Mulal81}. 
In this regard, let us 
mention that unlike in the latter reference, the upper cutoff 
at $w=w_u$ appears naturally in our formula.
It would also be interesting to
investigate the order of  contribution of more than two
 spinons to the
DCF, and in particular that of four spinons. Also of interest
is to find to what extent the two-spinon DCF can satisfy the
various sum rules  
which involve the full DCF \cite{Mulal81}. The extension of
this work to the Heisenberg  model with higher spin is certainly
desirable. In this case, the form factors can in 
principle be computed through the bosonization of the vertex
operators which is now available in \cite{Bou95}.
The full details of this paper and the relations with existing
approximate results and exact results of other models 
will appear elsewhere. 
We hope that 
the analysis
in this paper will make it convincing that infinite-dimensional
symmetries, and in particular infinite-dimensional quantum
groups might be of great usefulness in tackling some of the 
problems that have so far remained intractable and elusive 
 through 
traditional methods. Finally, we hope that this work will also
help improve interpretations of  experimental and
numerical data through the
spinon picture.

\newpage

\section{Aknowlegements}

The work of A.H.B. is supported by the NSF Grant \# PHY9309888.
We are  grateful to   
 Korepin and Weston for interesting discussions
at early stages of this work. We are particularly
thankful to Karbach and M\"uller for  
a stimulating discussion on the overall coefficient in 
formula (\ref{exact}). 
We highly appreciate the
help from Aurag, Kidonakis, Mateev and Smith  regarding  
numerical methods of integrations. 
We are also indebted to Klumper, McCoy, Perk,
 Sebbar, Shrock, 
Smirnov, and Takhtajan  for illuminating discussions.

\pagebreak

\end{document}